\documentclass[11pt]{article}

\usepackage{latexsym,color,amsmath,amsthm,amssymb,amscd,amsfonts}

\newtheorem{theo}{Theorem}

\newtheorem{prop}[theo]{Proposition}

\newtheorem{theorem}{Theorem}

\newtheorem{lemma}{Lemma}

\newtheorem{proposition}{Proposition}
\newtheorem{definition}{Definition}
\newtheorem{example}{Example}

\def\R{{\mathbb R}}
\def\K{{{\mathcal K}_n}}

\def\inte{{\rm int}}

\begin{document}

\title{Affine surface area
}

\author{Carsten Sch\"utt  and Elisabeth M. Werner 
\thanks{Partially supported by an NSF grant  DMS-2103482 }}

\date{}

\maketitle

\date{}

 \maketitle

\begin{abstract}
We give an overview of the affine surface area, its properties and its history.
\end{abstract}
\vskip 4mm

\section{Introduction}

A convex body in $\mathbb R^{n}$ is a convex, compact subset of $\mathbb R^{n}$
with nonempty interior. The set of all convex bodies in $\mathbb R^{n}$ is denoted by
$\mathcal K^{n}$. We denote the Euclidean norm of a vector $x\in \mathbb R^{n}$ by $\|x \|_{2}$.
The Euclidean ball in $\mathbb R^{n}$ with center $x$ and radius $r$ is $B_{2}^{n}(x,r)$.
We denote $B_{2}^{n}(0,r)$ by $B_{2}^{n}$.
We consider here $\mathcal K^{n}$ equipped with the Hausdorff metric
$$
d_{H}(C,K)
=\inf\{\rho| C\subseteq K+\rho B_{2}^{n}\hskip 1mm\mbox{and}\hskip 1mm
K\subseteq C+\rho B_{2}^{n}\}.
$$
Another metric on $\mathcal K^{n}$ is the symmetric difference metric
$$
d_{S}(C,K)=\operatorname{vol}_{n}(C\triangle K)
=\operatorname{vol}_{n}((C\setminus K)\cup(K\setminus C)).
$$
The support function $h_{K}:\mathbb R^{n}\to \mathbb R$ 
of a convex body $K$ in $\mathbb R^{n}$ is given by
$$
h_{K}(x)=\max_{y\in K}\langle x,y\rangle.
$$
The polar body of a convex body $C$ that contains the origin as an interior point is
$$
C^{*}=\{y\in\mathbb R^{n}|\forall x\in C:\hskip 1mm \langle x,y\rangle\leq1\}.
$$
The $n-1$-dimensional Hausdorff measure on $\mathbb R^{n}$ 
is denoted by $\mathcal H_{n-1}$ .
For a convex body $K$ in $\mathbb R^{n}$ the measure $\mu_{\partial K}$ is
 the restriction of $\mathcal H_{n-1}$ to the boundary of $K$. We call it
also the surface measure of $\partial K$.
$\sigma_{n}$ is the restriction of $\mathcal H_{n-1}$ to the Euclidean sphere $S^{n-1}$.
\par
The surface area measure 
\begin{equation}\label{SurfMeasure}
\sigma_{K} \hskip 2mm\mbox{on}\hskip 2mm S^{n-1}
\end{equation} 
is defined in the following way: For every Borel set
$A$ of $S^{n-1}$ we define
$\sigma_{K}(A)$ as the $n-1$-dimensional Hausdorff measure of the set
of all points in $\partial K$ that have a normal that is element in $A$.
\par
The boundary of a convex body or
a function mapping from $\mathbb R^{n}$ to $\mathbb R$ are of class $C^{2}$ if they are twice continuously differentiable
and they are of class $C^{2}_{+}$ is they are twice continuously differentiable and its Gau{\ss}-Kronecker
is strictly positive.
\par
For all $t_{1},\dots,t_{m}\geq0$ and all $K_{1},\dots,K_{m}\in\mathcal K^{n}$
(see e.g. Schneider \cite{SchneiderBuch})
$$
\operatorname{vol}_{n}(t_{1}K_{1}+\cdots+t_{m}K_{m})
=\sum_{i_{1},\dots,i_{n}=1}^{m}V(K_{i_{1}},\dots,K_{i_{n}})t_{i_{1}}\cdots t_{i_{n}}.
$$
The coefficients $V(\cdots)$ are called {\em mixed volumes}.
We denote
$$
V_{1}(K,C)=V(K,\dots,K,C)
$$
and we have
\begin{equation}\label{MixVol1}
V_{1}(K,C)
=\frac{1}{n}\lim_{\epsilon\to0}
\frac{\operatorname{vol}_{n}(K+\epsilon C)-\operatorname{vol}_{n}(K)}{\epsilon}.
\end{equation}
Moreover,
\begin{equation}\label{MixVol2}
V_{1}(K,C)=\frac{1}{n}\int_{S^{n-1}}h_{C}(\xi) d\sigma_{K}(\xi).
\end{equation}
We introduce the affine surface area and its main properties. We discuss
its history.
For some results we give the main steps of their proofs. For more background on general Convex Geometry we refer
to the books of Rockefellar \cite{Roc} and Schneider \cite{SchneiderBuch}.
All details of this paper can be found in the forthcoming book by Sch\"utt and Werner
\cite{SWBook}.

\section{Generalized Gau{\ss}-Kronecker curvature}

A convex body in $\mathbb R^{n}$ may not be twice differentiable at any boundary
point and therefore we do not have the classical notion of curvature. Thus we need to introduce a more general notion of second order differentiability.

\begin{definition}\index{subdifferential}
Let $\mathcal U$ be a convex, open subset of $\mathbb R^{n}$ and let 
$f:\mathcal U \rightarrow \mathbb R$ be a convex function. 
A vector $df(x_{0}) \in  \mathbb R^{n}$ is called 
subdifferential at the point $x_{0} \in \mathcal U$, if we have for all 
$x \in \mathcal U$
$$
f(x_{0})+\langle df(x_{0}),x-x_{0}\rangle \leq f(x).
$$
\end{definition}
\vskip 3mm

\begin{lemma}\label{SubdiffExist}
Let $\mathcal U$ be a convex, open set in $\Bbb R^{n}$ and let 
$f:\mathcal U \rightarrow \Bbb R$ be a convex function. Then $f$ has a subdifferential at
every $x \in \mathcal U$ and the set of subdifferentials at a given point $x$ is convex.
\end{lemma}

The existence of a subdifferential follows from the theorem of Hahn-Banach.

If a function $f$ is differentiable at a point $x$, then there is a unique subdifferential at $x$ and it is equal to its gradient.
By a theorem of Rademacher a convex function is  almost everywhere differentiable
(\cite{Roc}, p. 246), but there are convex functions 
that are nowhere twice differentiable. We give an example.

\begin{example}\label{ExDiff1}
Let $q:\mathbb N\to\mathbb Q\cap[0,1]$ be a bijection
 and let $g:[0,1] \rightarrow \Bbb R$ be given by
$$
g(x)=\sum_{q(n) \leq x}\frac{1}{2^{n}}.
$$
Moreover, let $G:[0,1]\to\mathbb R$ be the antiderivative of $g$
$$
G(x)=\int_{0}^{x}g(t)dt.
$$
Then $G$ is convex and differentiable at all irrational points, 
but not differentiable
at all rational points. The derivative of $G$ at the irrational point $x$
is $g(x_{})$.
In particular, $G$ is nowhere twice differentiable in the classical sense.
\end{example}

Since we need a notion of curvature for convex functions we introduce
the generalized second derivative of a convex function.

\begin{definition}\label{DefTwiceDiff1}
\index{function! twice differentiable}
Let $\mathcal U$ be a convex, open subset of $\mathbb R^{n}$
and $f:\mathcal U\to\mathbb R$ a convex function.
The convex function $f$ is said to be twice differentiable in the generalized
sense at $x_{0}$, if there are a linear map 
$d^{2}f(x_{0}):\mathbb R^{n}\to\mathbb R^{n}$ and a neighborhood
$\mathcal{V}(x_{0})$ such that we have for all $x \in \mathcal{V}(x_{0})$ and
for all subdifferentials $df(x)$
\begin{equation}\label{DefTwiceDiff1-1}
\left\|df(x)-df(x_{0})-d^{2}f(x_{0})(x-x_{0})\right\|_{2} \leq
\Theta(\|x-x_{0}\|_{2})\|x-x_{0}\|_{2},
\end{equation}
where $\Theta:[0,\infty)\to\mathbb R$ is a monotone function with $\lim_{t \to 0} \Theta(t)=0$.
The matrix of $d^{2}f(x_{0})$ with respect to the standard basis
of $\mathbb R^{n}$ is called generalized Hesse matrix.
Ignoring the difference between a linear map and its representation as a matrix we refer to $d^{2}f(x_{0})$ as the Hesse matrix.
\end{definition} 

Of course, if $f$ is twice differentiable in the classical sense then $d^{2}f(x_{0})$ equals the usual Hesse matrix
$\nabla^{2}f(x_{0})$. We give now an example of a function that is not twice differentiable, but that
is twice differentiable in the generalized sense.

\begin{example}\label{ExFuncSecondD1}
Let $f:[-1,1]\to\mathbb R$ be defined by
$$
f(x)=\left\{
\aligned
& x^{2} \hskip 35mm \text{if} \hskip 5mm |x|= \frac{1}{n}   \\
& \frac{2n+1}{n(n+1)}|x|-\frac{1}{n(n+1)}\hskip 5mm \text{if} \hskip 5mm
\frac{1}{n+1} < |x| < \frac{1}{n}
\endaligned
\right.
\hskip 10mm
n\in\mathbb N.
$$
Then $f$ is convex, twice differentiable at $0$ in the generalized sense and 
$d^{2}f(0)=2$, but $f$ is not twice differentiable at $0$ in the classical sense.
\end{example}

Most importantly, by a theorem of Busemann-Feller-Aleksandrov a convex function is almost everywhere
twice differentiable in the generalized sense \cite{Ale,Ban,BiCoPu,BuFe} and 
(\cite{GruberBuch}, p. 28).

We now define the {\em generalized Gau{\ss}-Kronecker curvature} of a convex function.
Let $\mathcal U$ be an open, convex set in $\mathbb R^{n}$ and 
$f:\mathcal U \rightarrow\mathbb R$ a convex function. 
The {\em epigraph} of $f$ is 
$$
\operatorname{epi}(f)=\{(x,t)|f(x)\leq t\}.
$$
For $x\in \mathcal U$ we put $z=(x,f(x))$. We assume that $\operatorname{epi}(f)$
has a unique normal $N(z_{0})$ at $z_{0}$.
Let 
$P_{N(z_{0})}:\mathbb R^{n+1}\to H(z_{0},N(z_{0}))$
be the orthogonal projection. If the 
convex sets 
\begin{equation}\label{DefGaussCurv-1}
\frac{1}{\sqrt{2\Delta}}P_{N(z_{0})}(\operatorname{epi}(f)\cap H(z_{0}-\Delta N(z_{0}),N(z_{0}))
\end{equation}
converge pointwise for $\Delta\to0$
to an ellipsoid or an elliptic cylinder we call its
boundary the {\em indicatrix of Dupin} and denote it by
 $\operatorname{Dupin}(f,z_{0})$.
\par
We call the squares of the lengths of the principal radii $r_{1},\dots,r_{n}$  of the indicatrix
of Dupin the {\em generalized principal Gau{\ss}-Kronecker curvature radii}
and, if $r_{1},\dots,r_{n}>0$, their reciprocals the {\em generalized principal Gau{\ss}-Kronecker curvatures}.
Then we say that the {\em generalized Gau{\ss}-Kronecker curvature} of $f$ at
$x_{0}$ exists and it is equal to the product  of the principal curvatures 
\begin{equation}\label{AffSurfPrincAx1}
\kappa(x_{0})=\left(\prod_{i=1}^{n}r_{i}\right)^{-2}
=\left(\frac{\operatorname{vol}_{n}(B_{2}^{n})}{\operatorname{vol}_{n}(\operatorname{Dupin}(f,z_{0}))}\right)^{-2}.
\end{equation}
When the indicatrix of Dupin is an elliptic cylinder we define the curvature to be $0$.
\par
One can show that all eigenvalues of the generalized Hesse matrix $d^{2}f(x_{0})$ are nonnegative.
If all eigenvalues are strictly positive, the indicatrix of Dupin is an elliptic sphere. In this case the sets \eqref{DefGaussCurv-1} converge uniformly in the Hausdorff metric to the elliptic sphere $\operatorname{Dupin}(f,z_{0}))$ because all sets are convex and bounded.
\par
It can be shown that the generalized Gau{\ss}-Kronecker curvature of a convex function
at a point $x_{0}$ at which $f$ is twice differentiable in the generalized sense equals
\begin{equation}\label{curvature1}
\kappa(x_{0})
=
\frac{\det(d^{2}f(x_{0}))}{(1+\|df(x_{0})\|_{2}^{2})^{\frac{n+2}{2}}}.
\end{equation}
In order to define the curvature of a convex body $K$ at $x_{0}\in\partial K$, we parametrize the boundary in a
neighborhood of $x_{0}$ by a convex function and take its curvature as the curvature of the convex body.

\section {Affine surface area}

\begin{definition}\index{affine surface area}
The affine surface area of a convex body $K$ in $\mathbb R^{n}$
is
\begin{equation}\label{DefAffSurf}
as(K)=\int_{\partial K}\kappa(x)^{\frac{1}{n+1}}d\mu_{\partial K}(x)
\end{equation}
where $\kappa$ denotes the generalized Gau{\ss}-Kronecker curvature \eqref{AffSurfPrincAx1}
and $\mu_{\partial K}$ the surface measure on $\partial K$.
\end{definition}

As noted above
by a theorem of Busemann-Feller-Aleksandrov the generalized curvature exists almost everywhere 
\cite{Ale,Ban,BiCoPu,BuFe} and thus
the integrand $\kappa^{\frac{1}{n+1}}$ is well defined, but we still have to confirm that
it is an integrable function. 

\par
We discuss the history of the affine surface area \cite{Leichtw4,Leichtw1,Lutw1}
and what led to the definition of the affine surface area \eqref{DefAffSurf} and its alternative definitions. In
Affine Differential Geometry the affine surface area of a convex body whose boundary has a $C^{2}_{+}$-parametrization $x:U\to\mathbb R^{n}$ was defined as (\cite{Bla2}, Chapter 47)
\begin{equation}\label{DefAaffSBla}
\int_{U}\left|\det\left(\det\left(\frac{\partial x}{\partial t_{1}},
\dots,\frac{\partial x}{\partial t_{n-1}},\frac{\partial^{2}x}{\partial t_{i}
\partial t_{j}}\right)_{i,j=1}^{n-1}\right)\right|^{\frac{1}{n+1}}d(t_{1},
\dots,t_{n-1}).
\end{equation}
This is the simplest expression that is independent of the change of parametrization,
is invariant under affine maps, that depends only on the surface of a convex body
and involves only first and second derivatives of the parametrization. One can show that
\eqref{DefAaffSBla} equals \eqref{DefAffSurf}.
\par
Thus \eqref{DefAaffSBla} is a definition of the affine surface area for all convex bodies with a $C^{2}_{+}$-boundary.
As the affine surface area has remarkable and useful properties - see below - it was imperative to extend it to arbitrary convex bodies.
With the properties of the affine surface area for convex bodies with $C_{+}^{2}$-boundary
in mind
we expect that the extended affine surface area satisfies
the following:
\vskip 1mm
\noindent
{\bf (i)} The extension should coincide with \eqref{DefAffSurf} and \eqref{DefAaffSBla} for convex bodies with $C^{2}_{+}$-boundary.
\vskip 1mm
\noindent
{\bf (ii)} For all polytopes $P$
$$
\operatorname{as}(P)=0.
$$
\vskip 1mm
\noindent
{\bf (iii)} The affine surface area is not continuous with respect to the Hausdorff distance. 
This follows from considering a sequence of polytopes that converge in the Hausdorff metric to the Euclidean ball $B_{2}^{n}$. But one wants to have upper semicontinuity.
\vskip 1mm
\noindent
{\bf (iv)} For all convex bodies $K$ in $\mathbb R^{n}$ and all affine maps $T:\mathbb R^{n}\to\mathbb R^{n}$
$$
\operatorname{as}(T(K))=|\det T|^{\frac{n-1}{n+1}}\operatorname{as}(K).
$$
\vskip 1mm
\noindent
{\bf (v)} The affine isoperimetric inequality should hold for all convex bodies $K$
in
$\mathbb R^{n}$
$$
\operatorname{as}_{}(K)^{n+1}
\leq \operatorname{as}_{}(B_{2}^{n})^{2}n^{n-1}\operatorname{vol}_{n}(K)^{n-1}.
$$
Since $\operatorname{as}_{}(B_{2}^{n})=\operatorname{vol}_{n-1}(\partial B_{2}^{n})$
this means in particular that $\operatorname{as}_{}(K)$
is smaller than or equal to the affine surface area of a Euclidean ball that has the same volume.
\vskip 1mm
\noindent
{\bf (vi)} The affine surface area appears naturally in approximation of 
convex bodies by polytopes.
That should also hold for the extended affine surface area.
\par

In order to give the affine surface area a geometric interpretation Blaschke introduced  the floating body (\cite{Bla2}, p. 125-6). 
He introduced the {\em floating body} of a convex body $K$ as this convex body $[K]_{t}$
whose tangent hyperplanes cut off a set of  exactly volume $t$ from
$K$. He assumed that $K$ has an analytic boundary to ensure that
$[K]_{t}$ exists, at least for small $t>0$. Then he showed in $\mathbb R^{3}$ 
\begin{equation}\label{FloatBod1}
\operatorname{as}(K)
=\lim_{t\to0}\frac{\operatorname{vol}_{3}(K)
-\operatorname{vol}_{3}([K]_{t})}{\sqrt{t}}.
\end{equation}
Leichtweiss generalized this to higher dimensions and the differentiability class
$C^{2}_{+}$ \cite{Leichtw1},
\begin{equation}\label{FloatBod11}
\operatorname{as}(K)
=\lim_{t\to0}\frac{\operatorname{vol}_{n}(K)
-\operatorname{vol}_{3}([K]_{t})}{t^{\frac{2}{n+1}}}.
\end{equation}
The name affine surface area reflects the Minkowski definition of the ususal surface area of a 
convex body
$$
\operatorname{vol}_{n-1}(\partial K)
=\lim_{t\to0}\frac{\operatorname{vol}_{n}( K+tB_{2}^{n})-\operatorname{vol}_{n}( K)}{t}.
$$
\par
Expression \eqref{FloatBod11}
 is the starting point for Leichtweiss for his extension of the affine surface area to arbitrary
convex bodies: Why not use the right hand side of \eqref{FloatBod11} to define the extension?
Unfortunately, the floating body may not exist in general. The simplex in $\mathbb R^{n}$ is an example.
Of course, for the limit to exist in \eqref{FloatBod11}, the floating body only needs to exist
for small $t\geq0$. Again, the simplex is a counterexample.
\par

In order to circumvent this difficulty, Leichtweiss \cite{Leichtw4} added a small Euclidean ball (or ellipsoid)
to the convex body. Now the floating body exists for small $t$ and therefore the affine surface area.
Leichtweiss then took the limit of the radius of the Euclidean ball to $0$
and defined this limit to be
 the affine surface area of the convex body $K$. Formally \cite{Leichtw4},
\begin{eqnarray}
&&n\cdot c_{n}\lim_{\epsilon\to0}\lim_{\delta\to0}\frac{1}{\delta^{\frac{2}{n+1}}}
(\operatorname{vol}_{n}(K+\epsilon\mathcal  E)-V(K+\epsilon\mathcal  E,\dots,K+\epsilon\mathcal  E,[K+\epsilon\mathcal  E]_{\delta}))
\nonumber\\
&&=n\cdot c_{n}\lim_{\epsilon\to0}\lim_{\delta\to0}\frac{1}{\delta^{\frac{2}{n+1}}}
(\operatorname{vol}_{n}(K+\epsilon\mathcal  E)-V_{1}(K+\epsilon\mathcal  E,
[K+\epsilon\mathcal  E]_{\delta})),
\label{DefAffLeicht}
\end{eqnarray}
where $[K+\epsilon\mathcal E]_{\delta}$ is the floating body of $K+\epsilon\mathcal E$.
$\mathcal E$ is an ellipsoid with $\operatorname{vol}_{n}(\mathcal E)=\operatorname{vol}_{n}(B_{2}^{n})$ and 
$$
c_{n}=2\left(\frac{\operatorname{vol}_{n}(B_{2}^{n})}{n+1}\right)^{\frac{2}{n+1}}.
$$
He showed that it is an extension of the affine surface area introduced by
Blaschke.
A somewhat simpler expression is (\cite{Leichtw5}, equation (11))
\begin{equation}\label{DefAffSA51}
c_{n}\lim_{\epsilon\to0}\lim_{\delta\to0}\frac{1}{\delta^{\frac{2}{n+1}}}
(\operatorname{vol}_{n}(K+\epsilon\mathcal  E)
-\operatorname{vol}_{n}([K+\epsilon\mathcal  E]_{\delta}))
\end{equation}
It was shown by Sch\"utt \cite{Schue3} that  \eqref{DefAffSurf} and \eqref{DefAffSA51} are equal.
\par
At the same time Lutwak gave another extension of the affine surface area  to arbitrary convex bodies by quite a different approach. A compact subset $L$ of $\mathbb R^{n}$ is star shaped, if
for all $x\in L$ the line segment $[0,x]$ is also contained in $L$. The radial function
$\rho_{L}:S^{n-1}\to\mathbb R$ of a star shaped set $L$ is defined by
$$
\rho_{L}(\xi)=\max\{\lambda\geq0|\lambda\xi\in L\}.
$$
If $\rho_{L}$ is continuous with respect to the Euclidean norm, we call $L$ a star body. $S_{c}^{n}$ denotes the set of all star bodies
$L$
whose centroid 
$$
\frac{1}{(n+1)\operatorname{vol}_{n}(L)}
\int_{S^{n-1}}\rho_{L}(\xi)^{n+1}\xi d\sigma_{n}(\xi)
$$
is at the origin.
In \cite{Lutw1} Lutwak introduced
\begin{equation}\label{DefLutwakAff}
n^{\frac{1}{n+1}}\inf_{L\in S^{n}_{c}}\left((\operatorname{vol}_{n}(L))^{\frac{1}{n}}
\int_{S^{n-1}}\frac{1}{\rho_{L}(\xi)}
d\sigma_{K}(\xi)\right)^{\frac{n}{n+1}}
\end{equation}
as the affine surface area.
The measure $\sigma_{K}$ 
on $\partial B_{2}^{n}$ is the surface area measure \eqref{SurfMeasure}.
\par
This definition is inspired by Petty's definition of the geominimal surface area
of a convex body \cite{Petty}. Petty defined the geominimal surface area as
\begin{equation}\label{GMSPetty}
\operatorname{gma}(K)
=\frac{n}{\operatorname{vol}_{n}(B_{2}^{n})^{\frac{1}{n}}}
\inf\left\{\left.V_{1}(K,C^{*})\operatorname{vol}_{n}(C)^{\frac{1}{n}}\right|C\hskip 1mm
\mbox{convex body with centoid at}\hskip 1mm 0\right\},
\end{equation}
where $C^{*}$ is the polar body to $C$ and $V_{1}(K,C^{*})$ is the mixed volume \eqref{MixVol1}.
\par
Again, at the same time Sch\"utt and Werner \cite{ScWe1} 
used the convex floating body
$K_{t}$ of a convex body $K$ to extend the affine surface area to arbitrary
convex bodies. Let $t\geq0$. The {\em convex floating body} $K_{t}$ of a convex body $K$ in $\mathbb R^{n}$
is the intersection of all halfspaces $H^{+}$ whose defining hyperplanes cut off a set of volume
$t$, i.e.
\begin{equation}\label{DefConvexFloat1}
K_{t}=\bigcap_{\operatorname{vol}_{n}(K\cap H^{-})=t} H^{+}.
\end{equation}
This notion has been independently introduced by B\'ar\'any and Larman \cite{BL}
and Sch\"utt and Werner \cite{ScWe1}.
\par
Sch\"utt and Werner showed for all convex bodies $K$ in $\mathbb R^{n}$ \cite{ScWe1}
\begin{equation}\label{AffSurfAr1}
\lim_{t \to 0} \frac{\operatorname{vol}_{n}(K)-\operatorname{vol}_{n}(K_{t})}{t^{\frac{2}{n+1}}}
=
\frac{1}{2}\left(\frac{n+1}{\operatorname{vol}_{n-1}(B_{2}^{n-1})}\right)^{\frac{2}{n+1}}
\int_{\partial K}\kappa(x)^{\frac{1}{n+1}}d\mu_{\partial K}(x),
\end{equation}
where $\kappa$ is the generalized Gau{\ss}-Kronecker curvature. 
\eqref{AffSurfAr1} suggests that we can use \eqref{DefAffSurf} as a definition for the
extended affine surface area. In fact, it can be shown that this expression satisfies
all the desired requirements.
\par
Werner proved that an analogous formula holds for the illumination body \cite{Werner1994}.

\section{Convex floating body and the rolling theorem}

Blaschke said that {\em a Euclidean ball rolls freely in a convex body $K$} if it can be placed
at each boundary point $x$ of $K$ such that it touches the boundary at $x$ and is contained 
in $K$. Formally, a Euclidean ball with radius $r$ rolls freely in $K$ if for all $x\in\partial K$
we have 
$$
B_{2}^{n}(x-rN(x))\subseteq K.
$$
If there is $r>0$ such that $rB_{2}^{n}$ rolls freely in $K$, then $K$ is of class $C^{1}$. Furthermore,
if there is $r>0$ such that $rB_{2}^{n}$ rolls freely in $K$, then $K$ is the $r$-outer parallel
body of another convex body $C$, i.e.
$$
K=C+rB_{2}^{n}.
$$
Note that there are convex bodies with a $C^{1}$-boundary such that no ball rolls freely.
The convex sets 
$$
B_{p}^{n}=\left\{x\in\mathbb R^{n}\left|\sum_{i=1}^{n}|x_{i}|^{p}\leq1\right.\right\}
$$ 
 for $1<p<2$ are examples. Furthermore,
for any convex body $K$ with a $C^{2}_{+}$-boundary there is a ball that rolls freely
in $K$, but not all bodies for which there is a ball that rolls freely have a
$C^{2}$ boundary. $C^{n}+B_{2}^{n}$ where $C^{n}$ is the $n$-dimensional cube
with sidelength $1$ is an example.

A quantitative version of {\em a ball rolls freely in the convex body} was
introduced by Sch\"utt and Werner \cite{ScWe1,SWBook}.
For $x \in \partial K$ we denote by $r(x)$ the supremum of all radii of
Euclidean balls that contain $x$ and that are contained in $K$, i.e.
$r:\partial K\to\mathbb R$ is defined by 
\begin{equation}\label{DefRollRad1}
r(x)=\sup\left\{\|x-z\||\exists z\in K:B_{2}^{n}(z,\|x-z\|)\subseteq K\right\}.
\end{equation}
We call $r$ the rolling function.\index{rolling function}
\index{function! rolling}
\par
If $K$ does not have a unique normal at $x$ then $r(x)=0$. Moreover, we have
\begin{equation}\label{DefRollRad2}
r(x)=\sup\{\rho|B_{2}^{n}(x-\rho N(x),\rho)\subseteq K\},
\end{equation}
if $K$ has a unique normal at $x$.

\begin{lemma}\cite{SWBook}
Let $K$ be a convex body in $\mathbb R^{n}$ and let 
$r:\partial K\to\mathbb R$ be the rolling function (\ref{DefRollRad1}).
Then we have 
\newline
(i) The supremum in  (\ref{DefRollRad1}) is attained.
\newline
(ii) $r$ 
is upper semicontinuous.
\end{lemma}

Howard \cite{Howard1} studied the rolling function on manifolds.
There are convex bodies for which the function $r$ is not continuous.
\vskip 3mm

\begin{example}
Let 
$$
K=\operatorname{conv}\left.\left\{\left(\pm\frac{1}{n},\frac{1}{n^{2}}\right)\right|n\in\mathbb N\right\}.
$$
and let $r$ be its rolling function.
Then $r((0,0))\geq\frac{1}{2}$ but for all $n\in\mathbb N$ we have
$r\left(\left(\frac{1}{n},\frac{1}{n^{2}}\right)\right)=0$. In particular $r$ is not continuous at
$(0,0)$.
\end{example}

Sch\"utt and Werner gave a quantitative result for {\em a ball rolls freely in a convex body}.

\begin{theorem}\label{SW1}  \cite{ScWe1} 
Let $K$ be a convex body in $\mathbb R^{n}$ 
such that it contains $B_{2}^{n}$. Then we have for all $t$ with $0 \leq t \leq 1$,
that $\{x \in \partial K|r(x) \geq t\}$ is a closed set and
$$
(1-t)^{n-1}\operatorname{vol}_{n-1}(\partial K) \leq
\mathcal H_{n-1}(\{x \in \partial K|r(x) \geq t\}).
$$
The inequality is optimal.
\par
In particular, the function $r^{-\alpha}:\partial K\to\mathbb R$ is Lebesgue integrable
for all $\alpha$ with $0\leq\alpha<1$.
\end{theorem}

Since $\kappa(x)^{\frac{1}{n-1}}\leq \frac{1}{r(x)}$ we have 
$\kappa(x)^{\frac{1}{n+1}}\leq r(x)^{-\frac{n-1}{n+1}}$.
By Theorem \ref{SW1}, the function $r^{-\frac{n-1}{n+1}}$ is integrable and consequently
$\kappa^{\frac{1}{n+1}}$ is smaller than an integrable function. This is used to prove
\eqref{AffSurfAr1}. Furthermore,
in order to prove \eqref{AffSurfAr1} we apply Lebesgue's dominated convergence theorem
and use $r^{-\frac{1}{n+1}}$ as the dominating function. 
\vskip 3mm

\noindent
{\bf Proof of \eqref{AffSurfAr1}.}
Assume that the orgin is an interior point of $K$. Then
$$
\operatorname{vol}_{n}(K)-\operatorname{vol}_{n}(K_{t})
=\frac{1}{n}\int_{\partial K}\langle x,N(x)\rangle
\left(1-\left|\frac{\|x_{t}\|_{2}}{\|x\|_{2}}\right|^{n}
\right)d\mu_{\partial K}(x),
$$
where $x\in\partial K$ and $x_{t}\in\partial K_{t}$ is the unique point in the intersection of the line segment
$[0,x]$ and $\partial K_{t}$. After dividing by $t^{\frac{2}{n+1}}$
$$
\frac{\operatorname{vol}_{n}(K)-\operatorname{vol}_{n}(K_{t})}{t^{\frac{2}{n+1}}}
=\frac{1}{n}\int_{\partial K}\frac{\langle x,N(x)\rangle}{t^{\frac{2}{n+1}}}
\left(1-\left|\frac{\|x_{t}\|_{2}}{\|x\|_{2}}\right|^{n}
\right)d\mu_{\partial K}(x).
$$
Now it is left to observe
$$
\kappa(x)^{\frac{1}{n+1}}
=
\lim_{t\to0}
\frac{\langle x,N(x)\rangle}{t^{\frac{2}{n+1}}}
\left(1-\left|\frac{\|x_{t}\|_{2}}{\|x\|_{2}}\right|^{n}
\right)
$$
and apply Lebesgue's dominated convergence theorem where we use
$r^{-\frac{n-1}{n+1}}$ as the dominating function.
$\Box$

\section{Properties of the affine surface area}

\begin{theorem}
The expressions \eqref{DefAffSurf}, \eqref{DefAffLeicht} and \eqref{DefLutwakAff} all coincide.
\end{theorem}

It was shown by Dolzmann and Hug
that the definitions of the affine surface area of Leichtweiss and Lutwak coincide   \cite{DolzHug}.

The projection body $\Pi(K)$ of a convex body $K$ is the convex body whose support function
is given by 
$$
h_{\Pi(K)}(\xi)=\operatorname{vol}_{n-1}(P_{\xi}(K)),
$$
where $P_{\xi}:\mathbb R^{n}\to H(0,\xi)$ is the orthogonal projection onto the hyperplane
containing $0$ and orthogonal to $\xi$.

\begin{theorem}\label{AffSurfAff}
(i) The affine surface area is an affine invariant, i.e.
for all affine maps $T:\mathbb R^{n}\to\mathbb R^{n}$ and all
convex bodies $K$ in $\mathbb R^{n}$
$$
\operatorname{as}( T(K))
=|\det(T)|^{\frac{n-1}{n+1}}\operatorname{as}(K).
$$
(ii) For all polytopes $P$ in $\mathbb R^{n}$ we have $\operatorname{as}(P)=0$.
\newline
(iii) The affine surface area $\operatorname{as}:\mathcal K^{n}\to\mathbb R$
is upper semicontinuous with respect to the Hausdorff metric on $\mathcal K^{n}$.
\newline
(iv) For all convex bodies $K$ in $\mathbb R^{n}$
\begin{equation}\label{AffSurfIneq}
\operatorname{as}(K)^{n+1}
\leq \operatorname{as}_{}(B_{2}^{n})^{2}n^{n-1}\operatorname{vol}_{n}(K)^{n-1}.
\end{equation}
Equality holds if and only if $K$ is an ellipsoid.
The inequality is called affine isoperimetric inequality.
\newline
(v) For all convex bodies $K$ in $\mathbb R^{n}$
$$
\operatorname{as}(K)^{n+1}
\leq \frac{n^{n+1}\operatorname{vol}_{n}(B_{2}^{n})^{n}}{\operatorname{vol}_{n-1}(B_{2}^{n-1})^{n}}
\operatorname{vol}_{n}(\Pi(K))
$$
with equality if and only if $K$ is an ellipsoid.
This inequality is called the Petty affine projection inequality.
\end{theorem}
\vskip 3mm

\noindent
{\bf Proof.}
(i) By \eqref{AffSurfAr1}
$$
\lim_{t \to 0} \frac{\operatorname{vol}_{n}(T(K))-\operatorname{vol}_{n}((T(K))_{t})}{t^{\frac{2}{n+1}}}
=
\frac{1}{2}\left(\frac{n+1}{\operatorname{vol}_{n-1}(B_{2}^{n-1})}\right)^{\frac{2}{n+1}}
\int_{\partial T(K)}\kappa(x)^{\frac{1}{n+1}}d\mu_{\partial T(K)}(x).
$$
Since $(T(K))_{t}=T(K_{\frac{t}{|\det(T)|}})$
\begin{eqnarray*}
&&
\frac{1}{2}\left(\frac{n+1}{\operatorname{vol}_{n-1}(B_{2}^{n-1})}\right)^{\frac{2}{n+1}}
\int_{\partial T(K)}\kappa(x)^{\frac{1}{n+1}}d\mu_{\partial T(K)}(x)
\\
&&=\lim_{t \to 0} \frac{\operatorname{vol}_{n}(T(K))-\operatorname{vol}_{n}((T(K_{\frac{t}{|\det(T)|}}))}{t^{\frac{2}{n+1}}}
=|\det(T)|\lim_{t \to 0} \frac{\operatorname{vol}_{n}(K)-\operatorname{vol}_{n}((K_{\frac{t}{|\det(T)|}})}{t^{\frac{2}{n+1}}}.
\end{eqnarray*}
With $s=\frac{t}{|\det(T)|}$
\begin{eqnarray*}
&&
\frac{1}{2}\left(\frac{n+1}{\operatorname{vol}_{n-1}(B_{2}^{n-1})}\right)^{\frac{2}{n+1}}
\int_{\partial T(K)}\kappa(x)^{\frac{1}{n+1}}d\mu_{\partial T(K)}(x)
=|\det(T)|^{\frac{n-1}{n+1}}\lim_{t \to 0} \frac{\operatorname{vol}_{n}(K)-\operatorname{vol}_{n}(K_{s})}{s^{\frac{2}{n+1}}}
\\
&&
=\frac{1}{2}|\det(T)|^{\frac{n-1}{n+1}}\left(\frac{n+1}{\operatorname{vol}_{n-1}(B_{2}^{n-1})}\right)^{\frac{2}{n+1}}
\int_{\partial K}\kappa(x)^{\frac{1}{n+1}}d\mu_{\partial K}(x).
\end{eqnarray*}
\par
(ii)
Since the curvature of a polytope is $0$ almost everywhere the equality follows 
immediately by \eqref{AffSurfAr1}.
Lutwak also showed this for his definition of affine surface area \eqref{DefLutwakAff}
in \cite{Lutw1}.
\par
(iii) Lutwak showed that the affine surface area is upper semicontinuous \cite{Lutw1}.
Another proof of the upper semicontinuity by Ludwig can be found in \cite{Ludw1}.
\par
(iv) Lutwak proved this inequality (\cite{Lutw1}, Corollary 7.7). 
\par
Hug gave a proof using Steiner symmetrization \cite{Hu1}. It is very similar to the proof
of the isoperimetric inequality via Steiner symmetrization.
The same 
idea appears in the proof of the Blaschke-Santal\'o inequality by Meyer and Pajor \cite{MePa}.
\par
Let $K$ be a convex body in $\mathbb R^{n}$ that contains $0$ as an interior point and let $H(0,\xi)$ be the hyperplane containing $0$ and being orthogonal to $\xi$.
Let $P_{\xi}$ be the orthogonal projection onto the hyperplane $H$. We define 
$f_{\xi}^{-}, f_{\xi}^{+}:P_{\xi}(K)\to\mathbb R$ by
\begin{eqnarray*}
f_{\xi}^{-}&=&\min\{t\in\mathbb R|x+t\xi\in K\}
\label{Param1}
\\
f_{\xi}^{+}&=&\max\{t\in\mathbb R|x+t\xi\in K\}
\label{Param2}
\end{eqnarray*}
By \eqref{curvature1}
$$
\operatorname{as}(K)
=\int_{P_{\xi}(K)}(\det(d^{2}f_{\xi}^{-}))^{\frac{1}{n+1}}
d_{H}x
+\int_{P_{\xi}(K)}(\det(d^{2}(-f_{\xi}^{+})))^{\frac{1}{n+1}}d_{H}x,
$$
where $d_{H}x$ is the Lebesgue measure on the hyperplane $H$.
For all symmetric, positive semidefinite
$(n-1)\times(n-1)$-matrices $A$ and $B$
$$
(\det A)^{\frac{1}{n+1}}+(\det B)^{\frac{1}{n+1}}
\leq2\left(\det \left(\frac{1}{2}(A+B)\right)\right)^{\frac{1}{n+1}}.
$$
This inequality follows from the standard inequality for positive semidefinite $n\times n$-matrices
$A$ and $B$
$$
(\det A)^{\frac{1}{n}}+(\det B)^{\frac{1}{n}}
\leq\left(\det (A+B)\right)^{\frac{1}{n}}
$$
by considering the $(n+1)\times(n+1)$-matrices
$$
\left(\begin{array}{ccc}
A & &  \\
 & 1 &  \\
 &  & 1
\end{array}\right)
\hskip 10mm\mbox{and}\hskip 10mm
\left(\begin{array}{ccc}
B& &  \\
 & 1 &  \\
 &  & 1
\end{array}\right).
$$
Therefore
$$
\operatorname{as}(K)
\leq2\int_{p_{H}(K)}\left(\det\left(\frac{1}{2}d^{2}(f_{\xi}^{-}-f_{\xi}^{+})\right)\right)^{\frac{1}{n+1}}
d_{H}x
=\operatorname{as}(\operatorname{St}_{H}(K)).
$$
Now we choose a sequence of symmetrizations that transforms the convex body
$K$ into the Euclidean ball. Finally, we use that the affine surface area is upper
semicontinuous.
\par
(v) Lutwak proved the Petty affine projection inequality for arbitrary convex bodies \cite{Lutw1}.
$\Box$
\vskip 3mm

In general, it is not easy to compute the affine surface area of a given convex body.
But it is possible for the convex bodies $B_{p}^{n}$.
For this we need the following lemma.

Let $A$ be a $n\times n$-matrix. The $(k,\ell)$-th coordinate of the cofactor matrix \index{matrix! cofactor}
\index{cofactor matrix}
$\operatorname{cof}(A)$ is $(-1)^{k+\ell}$ times
the determinant
of the matrix $A$ after deleting the $k$-th row and 
$\ell$-th column.

\begin{lemma}\cite{Domb,Gold, Gromoll}
\label{CurvImplFunc1}
Let $\mathcal U$ be a convex, open subset of $\mathbb R^{n+1}$ and $x_{0}\in\mathcal U$.
Let $F:\mathcal U\to\mathbb R$ be
twice continuously differentiable in a neighborhood of $x_{0}$. Moreover suppose that $F(x_{0})=0$ and that $F(x)\leq0$ is a convex body contained in $\mathcal U$. Then
\begin{equation}\label{CurvImplFunc1-1}
\kappa(x_{0})=\frac{|(\nabla F(x_{0}))^{t}\operatorname{cof}(\nabla^{2} F)(x_{0}) \nabla F(x_{0})|}
{\|\nabla F(x_{0})\|_{2}^{n+2}}
=\frac{\left|\det\left(
\begin{array}{cc}
\nabla^{2}F & \nabla F
\\
(\nabla F)^{t}  &0
\end{array}\right)\right|}{\|\nabla F\|_{2}^{n+2}},
\end{equation}
where $\nabla F$ denotes the gradient and $\nabla^{2}F$ the Hessian of $F$.
\end{lemma}

\begin{example}\index{affine surface area! unit ball of $\ell_{p}^{n}$}
 Let $1<p<\infty$. Then we have for the Gau{\ss}-Kronecker curvature of $\partial B_{p}^{n}$
at $x$ with $x_{i}\ne0$ for all $i=1,\dots,n$ 
$$
\kappa(x)
=\frac{(p-1)^{n-1}(\prod_{i=1}^{n} |x_{i}|^{p-2})}
{(\sum_{i=1}^{n} |x_{i}|^{2p-2})^{\frac{n+1}{2}}}.
$$
And for the affine surface area we get
$$
\int_{\partial B_{p}^{n}} \kappa(x)^{\frac{1}{n+1}}d\mu_{\partial B_{p}^{n}}(x)
=\frac{2^{n}(p-1)^{\frac{n-1}{n+1}}(
\Gamma(\frac{p+n-1}{(n+1)p})^{n}}{p^{n-1}\Gamma(\frac{n(p+n-1)}{(n+1)p})}.
$$
\end{example}
\vskip 3mm

A map $\Phi:\mathcal K^{n}\to\mathbb R$ is called a valuation if
for all $K,C\in\mathcal K^{n}$ such that $K\cup C\in\mathcal K^{n}$ 
$$
\Phi(K\cup C)+\Phi(K\cap C)=\Phi(K)+\Phi(C).
$$

\begin{theorem}\cite{Schue3}\label{AffSurfVal1}
\index{affine surface area! valuation}
 The affine surface area is a valuation, i.e.
 for all $K,C\in\mathcal K^{n}$ such that $K\cup C\in\mathcal K^{n}$ 
 \begin{equation}\label{AffSurfVal1-1}
 \operatorname{as}(K\cup C)+ \operatorname{as}(K\cap C)
 = \operatorname{as}(K)+ \operatorname{as}(C).
 \end{equation}
\end{theorem}

\noindent
{\bf Proof.}
We decompose the involved boundaries into disjoint sets
\begin{eqnarray*}
\partial(C\cup K)
&=&(\partial C\cap\partial K)\cup(\partial C\cap K^{c})\cup(\partial K\cap C^{c})
\\
\partial(C\cap K)
&=&(\partial C\cap\partial K)\cup(\partial C\cap\operatorname{int}(K))\cup(\partial K\cap \operatorname{int}(C))
\\
\partial C&=&(\partial C\cap\partial K)\cup(\partial C\cap K^{c})\cup(\partial C\cap\operatorname{int}(K))
\\
\partial K&=&(\partial C\cap\partial K)\cup(\partial K\cap C^{c})\cup(\partial K\cap\operatorname{int}(C)).
\end{eqnarray*}
In order to show \eqref{AffSurfVal1-1} it is enough to prove
\begin{eqnarray*}
&&\int_{\partial C\cap\partial K}\kappa_{C\cup K}^{\frac{1}{n+1}}d\mu_{\partial (C\cup K)}
+\int_{\partial C\cap\partial K}\kappa_{C\cap K}^{\frac{1}{n+1}}d\mu_{\partial (C\cup K)}
\\
&&=\int_{\partial C\cap\partial K}\kappa_{C}^{\frac{1}{n+1}}d\mu_{\partial (C)}
+\int_{\partial C\cap\partial K}\kappa_{ K}^{\frac{1}{n+1}}d\mu_{\partial (C)},
\end{eqnarray*}
where $\kappa_{C\cup K}$ is the curvature of $\partial (C\cup K)$ and the others accordingly.
This follows from
$$
\kappa_{C\cup K}(x)=\min\{\kappa_{C}(x),\kappa_{K}(x)\}
\hskip 20mm
\kappa_{C\cap K}(x)=\max\{\kappa_{C}(x),\kappa_{K}(x)\}.
$$
This in turn is true since the indicatrix of Dupin \eqref{DefGaussCurv-1} of $C\cup K$ is the union of the
indicatrices of Dupin  of $C$ and $K$ at $x$ and the indicatrix of Dupin 
of $C\cap K$
at $x$ is the intersection of the indicatrices of $C$ and $K$ at $x$. Moreover,
the intersection or union of two ellipsoids is again an ellipsoid if and only if one ellipsoid is contained in the other.
$\Box$

By results of Hadwiger a map $\Phi:\mathcal K^{n}\to\mathbb R$ is a continuous valuations
that is invariant under rigid motions 
if and only if it is a linear combination of quermassintegrals. 
 A continuous valuation that is invariant under rigid motions and that is
$n-j$-homogeneous is a multiple of the $j$-th quermassintegral \cite{Hadwiger}.
\par
Ludwig and Reitzner \cite{LudReitz1} showed that aside from 
 the Euler characteristc and
the volume, there is only one nontrivial upper semi-continuous valuation
that is affinely invariant, namely
the affine surface area.

\begin{theorem} \cite{LudReitz1}\label{LudReit1}
A functional $\Phi:\mathcal K^{n}\to\mathbb R$ is an upper 
semicontinuous, translations and $SL(n)$-invariant valuation
if there are constants $c_{0},c_{1}\in\mathbb R$ and $c_{2}\geq0$
such that for all $K\in\mathcal K^{n}$
$$
\Phi(K)=c_{0}V_{0}(K)+c_{1}\operatorname{vol}_{n}(K)+c_{2}\operatorname{as}(K),
$$
where $V_{0}$ is the Euler characteristic.
\end{theorem}

\section{Random polytopes and best approximation}

The affine surface area appears in 
best approximation of convex bodies by polytopes and in approximation
by random polytopes. The next theorems are examples.

\begin{theorem}
Let $K$ be a convex body in $\mathbb R^{n}$ that has a $C_{+}^{2}$-boundary.
Then there are constants $\operatorname{del}_{n-1}$ and $\operatorname{div}_{n-1}$
depending only on the dimension $n$ such that
\begin{equation}\label{}
\lim_{N\to\infty}\left(N^{\frac{2}{n-1}}\inf_{\operatorname{vert}(P)\leq N\atop
P\subseteq K}d_{S}(K,P)\right)
=\frac{1}{2}\operatorname{del}_{n-1}
\left(\int_{\partial K}\kappa(x)^{\frac{1}{n+1}}d\mu_{\partial K}(x)
\right)^{\frac{n+1}{n-1}}
\end{equation}
\begin{equation}\label{}
\lim_{N\to\infty}\left(N^{\frac{2}{n-1}}\inf_{\operatorname{vert}(P)\leq N\atop
K\subseteq P}d_{S}(K,P)\right)
=\frac{1}{2}\operatorname{div}_{n-1}
\left(\int_{\partial K}\kappa(x)^{\frac{1}{n+1}}d\mu_{\partial K}(x)
\right)^{\frac{n+1}{n-1}},
\end{equation}
where $P$ denotes polytopes.
\end{theorem}

This theorem was proved by McClure and Vitale in dimension $2$ \cite{McVi}. Gruber proved this for arbitrary dimensions \cite{Grub2}.
\par
In order to estimate or compute the constants $\operatorname{del}_{n-1}$ it is enough
to consider the case $K=B_{2}^{n}$. Kabatjanskii and Levenstein \cite{KabLev} estimated
$\operatorname{del}_{n-1}$ from above, Gordon, Reisner and Sch\"utt from below
\cite{GRS1,GRS2}. The proofs were simplified and the estimates improved by
Mankiewicz and Sch\"utt  \cite{MaS1,MaS2}
$$
\frac{n-1}{n+1}\operatorname{vol}_{n-1}( B_{2}^{n-1})^{-\frac{2}{n-1}}
\leq
\operatorname{del}_{n-1}
\leq 2^{0.802}\operatorname{vol}_{n-1}(\partial B_{2}^{n})^{-\frac{2}{n-1}}
$$
and 
$$
\lim_{n\to\infty}\frac{\operatorname{del}_{n-1}}{n}=\frac{1}{2\pi e}
=0.0585498...
$$
\par
A random polytope of a convex body $K$ is the convex hull of finitely many
points that are chosen randomly from $K$ with respect to a probability measure.
The expected volume of a random polytope of $N$ chosen points is
\begin{equation}\label{RandPoly1}
\mathbb E(K,N)
=\left(\frac{1}{\operatorname{vol}_{n}(K)}\right)^{N}
\int_{K}\cdots\int_{K}\operatorname{vol}_{n}([x_{1},\dots,x_{N}])dx_{1}\cdots dx_{N},
\end{equation}
where $[x_{1},\dots,x_{N}]$ denotes the convex hull of the points
$x_{1},\dots,x_{N}$.

\begin{theorem}\label{RanPol1}
Let $K$ be a convex body in $\mathbb R^{n}$. Then
\begin{equation}
c_{n}\lim_{N\to\infty}\frac{\operatorname{vol}_{n}(K)
-\mathbb E(K,N)}{\left(\frac{\operatorname{vol}_{n}(K)}{N}\right)^{\frac{2}{n+1}}}
=\int_{\partial K}\kappa^{\frac{1}{n+1}}d\mu_{\partial K},
\end{equation}
where
$$
c_{n}=2\left(\frac{\operatorname{vol}_{n-1}(B_{2}^{n-1})}{n+1}\right)^{\frac{2}{n+1}}
\frac{(n+3)(n+1)!}{(n^{2}+n+2)(n^{2}+1)\Gamma\left(\frac{n^{2}+1}{n+1}\right)}.
$$
\end{theorem}

Theorem \ref{RanPol1} was first proved by B\'ar\'any for convex bodies with 
$C^{3}$-boundary \cite{Ba1}. For general convex bodies it was proved by Sch\"utt \cite{Schue1},
see also \cite{BHH1}.

The following theorem answers the analogous question for random polytopes whose vertices are chosen from the boundary of the convex body. Let $f:\partial K\to\mathbb R$ be a continuous
function with $\int_{\partial K}fd\mu_{\partial K}=1$. Then the expected volume
of a random polytope of $N$ points chosen randomly from the boundary of $K$ 
with respect to the probability measure $fd\mu_{\partial K}$ is
$$
\mathbb E(f,K)
=\int_{\partial K}\cdots\int_{\partial K}\operatorname{vol_{n}([x_{1},\dots,x_{N}])}
f(x_{1})\cdots f(x_{N})d\mu_{\partial K}\cdots d\mu_{\partial K}.
$$

\begin{theorem}\cite{SchueW4}\label{RanPol2}
Let $K$ be a convex body in $\mathbb R^{n}$ such that there are $r$ and $R$ in $\mathbb R$ with $0<r\leq R<\infty$ so that we have for all $x\in\partial K$
$$
B_{2}^{n}(x-rN(x),r)\subseteq K
\subseteq B_{2}^{n}(x-RN(x),R)
$$
and let $f:\partial K\to [0,\infty)$ be a continuous function
with $\int_{\partial K}fd\mu_{\partial K}=1$. Let $\mathbb P_{f}$ be the probability measure on $\partial K$ with
$d\mathbb P_{f}(x)=f(x)d\mu_{\partial K}(x)$. Then we have
\begin{equation}\label{RanPol2-1}
\lim_{N\to\infty}\frac{\operatorname{vol}_{n}(K)
-\mathbb E(f,N)}
{(\frac{1}{N})^{\frac{2}{n-1}}}
=c_{n}\int_{\partial K}\frac{\kappa(x)^{\frac{1}{n-1}}}
{f(x)^{\frac{2}{n-1}}}d\mu_{\partial K}(x),
\end{equation}
where $\kappa$ is the generalized Gau{\ss}-Kronecker
curvature and
$$
c_{n}
=\frac{(n-1)^{\frac{n+1}{n-1}}\Gamma(n+1+\frac{2}{n-1})}{2(n+1)!(\operatorname{vol}_{n-2}(\partial B_{2}^{n-1}))^{\frac{2}{n-1}}}.
$$
The minimum at the right hand side is attained for the normalized affine surface area measure with density
$$
f_{\operatorname{as}}(x)
=\frac{\kappa(x)^{\frac{1}{n+1}}}
{\int_{\partial K}\kappa(x)^{\frac{1}{n+1}}d\mu_{\partial K}(x).
}
$$
\end{theorem}

This result was proved by Sch\"utt and Werner \cite{SchueW4} and at the same time for
convex bodies with $C_{+}^{2}$-boundary by Reitzner \cite{Rei1}.
It is interesting that the minimum of \eqref{RanPol2-1} is attained for the affine surface measure.

\section{Constrained convex bodies with maximal affine surface area}

We introduce the analogue to John's theorem, when volume is replaced by affine surface area. 
In parallel to John's maximal volume ellipsoid, Giladi, Huang, Sch\"utt and Werner investigated  these convex bodies  contained in \( K\)  that have the largest affine surface areas.

The isotropic constant 
$L_K$ of $K$  is defined by 
\begin{equation} \label{LK}
n L_K^2= \min \left\{\left.\frac{1}{\operatorname{vol}_{n}(TK)^{1 + \frac{2}{n}}} \int_{a+TK} \|x \| ^2 dx \right| a \in \mathbb{R}^n, T \in GL(n) \right\}.
\end{equation}
\par
\noindent
\begin{theorem} \cite{GHSW}
  There is a  constant $c>0$ such that for all $n \in \mathbb{N}$ and all convex bodies 
$K\subseteq \R^n$,  
\begin{align*}
\frac{1}{n^{5/6}}  \  \left(\frac{c}{L_K}\right)^{\frac{2n}{n+1}}  \  n \operatorname{vol}_{n}(B^n_2)^ \frac{2}{n+1}
\operatorname{vol}_{n}(K)^\frac{n-1}{n+1}
\leq \sup_{C\subseteq K}  \operatorname{as}(C)
\leq  \   n \operatorname{vol}_{n}(B^n_2)^ \frac{2}{n+1}\operatorname{vol}_{n}(K)^\frac{n-1}{n+1},
\end{align*}
where we take the supremum over all convex bodies $C$ that are contained in $K$.
Equality holds in the right inequality iff $K$ is a centered ellipsoid.
\end{theorem}
\vskip 2mm

The theorem shows that  $\sup_{C\subseteq K}  \operatorname{as}(C)$ is proportional to 
 $\operatorname{vol}_{n}(K)^\frac{n-1}{n+1}$.
\vskip 3mm

\noindent
{\bf Proof.}
The right hande side inequality follows immediately from the affine isoperimetric inequality
\eqref{AffSurfIneq}.
\par
We use a thin shell estimate by Gu\'edon and E. Milman \cite{GueMil}, see also Paouris \cite{Paouris2006}, on concentration of volume: Given an isotropic random vector $X$ with log-concave
density in Euclidean space $\mathbb R^{n}$, we have for all $t\geq0$
$$
\mathbb P(\|X\|_{2}-\sqrt{n}|\geq t\sqrt{n})
\leq C\exp(-c\sqrt{n}\min\{t^{3},t\}).
$$
For the estimate from below we choose as the convex body
$$
C=K\cap B_{2}^{n}(0,c_{n}),
$$
where $c_{n}$ is appropriately chosen and is of the order $\sqrt{n}$.
By the concentration result, the volume of $K\cap B_{2}^{n}(0,c_{n})$
is of the same order as $K$. This also implies that the surface area of
$K\cap B_{2}^{n}(0,c_{n})$ is of the same arder as that of $K$. Furthermore,
this implies that the part of the boundary of $K\cap B_{2}^{n}(0,c_{n})$
which is also part of the boundary of $B_{2}^{n}(0,c_{n})$ is big.
The curvature on that part is easily computed.
$\Box$


\vskip 5mm

\begin{thebibliography}{~~}


\bibitem{Ale} 
{\sc A.D. Aleksandrov},
{\em Almost everywhere existence of the second differential
   of a convex function and some properties of convex surfaces
   connected with it},
 Uchenye Zapiski Leningrad Gos. Univ., Math. Ser.
6 (1939),  
 3-35.
 


\bibitem{Ban}
{\sc V. Bangert}
{\em Analytische Eigenschaften konvexer Funktionen auf
Riemannschen Mannigfaltigkeiten},
Journal f\"ur die Reine und Angewandte Mathematik
307 (1979), 309-324.



\bibitem{BL}
{\sc I. B\'ar\'any and D.G. Larman},
{\em Convex bodies, economic cap covering, random
polytopes},
Mathematika 35 (1988), 274-291.


\bibitem{Ba1} 
{\sc I. B\'ar\'any},
{\em Random polytopes in smooth convex bodies},
Mathematika 39 (1992), 81-92.


\bibitem{ Ba2}
{\sc I. B\'ar\'any},
{\em Approximation by random polytopes is almost optimal},
 II International Conference in "Stochastic Geometry, Convex Bodies and 
Empirical Measures"(Agrigento 1996) 
 Rend. Circ. Mat. Palermo (2) Suppl. No. 50 (1997), 43--50.
 
 


\bibitem{BiCoPu}
{\sc G. Bianchi, A. Colesanti, and C. Pucci},
{\em On the second differentiability of convex surfaces},
 Geometriae Dedicata
60 (1996),  39-48.






\bibitem{Bla1}
{\sc W. Blaschke}, {\em {\"U}ber affine Geometrie VII. Neue Extremeigenschaften von
Ellipse und Ellipsoid}, Leipz. Ber. 69, (1917), 306-318.

\bibitem{Bla2}
{\sc W. Blaschke},
{\em Differentialgeometrie II,
Affine Differentialgeometrie}, Berlin,
Springer-Verlag 1923.




\bibitem{Bla1}
{\sc W. Blaschke},
{\em Integralgeometrie 2: Zu Ergebnissen von M.W. Crofton},
 Bull. Math. Soc. Roumaine Sci. 37 (1935), 3-11.
 
 \bibitem{BF1}
 {\sc T. Bonnesen and W. Fenchel},
 {\em Theorie der konvexen K\"orper},
 1934, Springer-Verlag.



\bibitem{BHH1}
{\sc K.J. B\"or\"oczky, L.M. Hoffmann and D. Hug},
{\em Expectation of intrinsic volumes of random polytopes},
Periodica Mathematica Hungarica 57 (2008), 143-164.



\bibitem{BrIv}
{\sc E. M. Bronshteyn and L. D. Ivanov},
{\em The approximation of convex sets by polyhedra},
 Siberian Math. J. 16 (1975), 852-853.



 
 \bibitem{BuFe}
 {\sc H. Busemann and W. Feller},
 {\em Kr\"ummungseigenschaften konvexer Fl\"achen},
 Acta Mathematica 66 (1935), 1-47.
 
 \bibitem{DoCarmo}
 {\sc M.P. Do Carmo},
 {\em Differential Geometry of Curves and Surfaces},
 Prentice-Hall, 1976.
 

\bibitem{DolzHug}
{\sc G. Dolzmann and D. Hug},
{\em Equality of two representations of extended affine surface area},
Archiv der Mathematik 65 (1995), 352-356.

\bibitem{Domb}
{\sc P. Dombrowski},
{\em Kr\"ummungsgr\"o{\ss}en gleichungsdefinierter Untermannigfaltigkeiten Riemannscher Mannigfaltigkeiten},
Mathematische Nachrichten 38 (1968), 133-180.

\bibitem{Dupin}
{\sc C. Dupin},
{\em Application de g\'eometrie et de m\'echanique \'a la marine,
aux ponts et chauss\'ees}, Paris 1822.





 







\bibitem{GHSW}
{\sc O. Giladi, H. Huang, C. Sch\"utt and E.Werner},
{\em Constrained convex bodies with extremal affine surface areas},
Journal of Functional Analysis 279 (2020).

\bibitem{Gold}
{\sc R. Goldman},
{\em Curvature formulas for implicit curves and surfaces},
Computer Aided Geometric Design 22 (2005), 632-658.




\bibitem{GRS1}
{\sc Y. Gordon, S. Reisner and C. Sch\"utt} 
{\em  Umbrellas
and Polytopal Approximation of the Euclidean Ball},
 Journal of Approximation Theory  90 (1997), 9-22.


\bibitem{GRS2} 
{\sc Y. Gordon, S. Reisner and C. Sch\"utt}, 
{\em  Erratum}, Journal of Approximation Theory 95 (1998), 331.


\bibitem{Gromoll}
{\sc D. Gromoll, W. Klingenberg and W. Meyer},
{\em Riemannsche Geometrie im Gro{\ss}en},
Lecture Notes in Mathematics 55 (1975), Springer-Verlag.

\bibitem{Grub1}
{\sc P.M. Gruber},
{\em Volume approximation of convex bodies by inscribed
polytopes}, Mathematische Annalen 281 (1988), 292-245.


\bibitem{Grub2}
{\sc P.M. Gruber},
{\em Asymptotic estimates for best and stepwise approximation of
convex bodies II}, Forum Mathematicum 5 (1993), 521-538.


\bibitem{Grub3}
{\sc  P.M. Gruber},
{\em Comparisons of best and random approximation of convex bodies
by polytopes}, Rend. Circ. Mat. Palermo, II. Ser., Suppl. 50 (1997), 189-216.


\bibitem{Grub4}
{\sc P.M. Gruber},
{\em Approximation of convex bodies},
Convexity and its Applications, Birk-h\"auser
 P.M. Gruber and J.M. Wills
 1983, 131-162.


\bibitem{Grub5}
{\sc P.M. Gruber},
{\em Aspects of approximation of convex bodies},
 P.M. Gruber and J.M. Wills:
 Handbook of Convex Geometry,
 North-Holland, 1993, 319-345.
 
  \bibitem{GruberBuch}
 {\sc P.M. Gruber},
 {\em Convex and Discrete Geometry},
 Springer-Verlag 2007.




\bibitem{GueMil}
{\sc O. Gu\'edon and E. Milman},
{\em Interpolating thin-shell and sharp
large-deviation estimates for isotropic log-concave measures},
Geometric and Functinal Analysis 21 (2011), 1043-1068.

\bibitem{Gugg1}
{\sc H.W. Guggenheimer},
{\em Differential Geometry},
Dover Publications, New York 1977.



\bibitem{Hadwiger}
{\sc H. Hadwiger},
{\em Vorlesungen \"uber Inhalt, Oberfl\"ache und Isoperimetrie},
Springer-Verlag 1957.

\bibitem{Howard1}
{\sc R. Howard},
{\em Blaschke's rolling theorem for manifolds with boundary},
manuscripta mathematica 99 (1999), 471-483.



\bibitem{Hu1}
{\sc D. Hug},
{\em Contributions to Affine Surface Area},
 Manuscripta Mathematica 91 (1996), 283-301.
 
 


\bibitem{J}
{\sc F. John},
{\em Extremum problems with inequalities as subsidiary conditions},
in: R. Courant Anniversary Volume, 1948,
Interscience New York, 187-204.


\bibitem{KabLev}
{\sc G.A. Kabatjanskii and V.I. Levenstein},
{\em Bounds for packings on a sphere and in space},
Problems Inform. Transmission 14 (1978), 1-17.








\bibitem{Leichtw4}
{\sc K. Leichtweiss},
{\em Zur Affinoberfl\"ache konvexer K\"orper},
Manuscripta Mathematica 56 (1986), 429-464.


\bibitem{Leichtw1}
{\sc K. Leichtweiss},
{\em \"Uber eine Formel Blaschkes zur Affinoberfl\"ache},
Studia Sci. Math. Hung. 21 (1987), 453-474.

\bibitem{Leichtw2}
{\sc K. Leichtweiss},
{\em \"Uber einige Eigenschaften der Affinoberfl\"ache beliebiger konvexer
K\"orper},
Resultate der Mathematik 13 (1988), 255-282.




\bibitem{Leichtw3}
{\sc K. Leichtweiss},
{\em Bemerkungen zur Definition einer erweiterten
Affinoberfl\"ache von E. Lutwak},
Manuscripta Mathematica 65 (1989), 181-197.

\bibitem{Leichtw5}
{\sc K. Leichtweiss},
{\em On the history of the affine surface area for convex bodies},
Results in Mathematics 20 (1991), 650-656.



\bibitem{Lei}
{\sc K. Leichtweiss},
{\em Convexity and differential geometry},
 Handbook of Convex Geometry
volume  B,
 P.M. Gruber and J.M. Wills editors, (1993),
 1045-1080

\bibitem{Leichtw5}
{\sc K. Leichtweiss},
{\em Affine Geometry of Convex Bodies},
Johann Ambrosius Barth, Heidelberg, 1998.


\bibitem{Ludw1}
{\sc M. Ludwig},
{\em On the semicontinuity of curvature integrals},
Mathematische Nachrichten 227 (2001), 99-108.

\bibitem{LudReitz1}
{\sc M. Ludwig and M. Reitzner},
{\em A characterization of affine surface area},
Advances of Mathematics 147 (1999), 138-172.


\bibitem{Lutw1}
{\sc E. Lutwak},
{\em Extended affine surface area},
Advances in Mathematics 85 (1991), 39-68.


\bibitem{Lutw2}
{\sc E. Lutwak},
{\em The Brunn-Minkowski-Fiery  Theory II:
Affine and Geominimal Surface Areas},
 Advances in Mathematics 118 (1996), 244-194.


\bibitem{Lutw3}
{\sc E. Lutwak},
{\em On the semicontinuity of curvatures},
Comment. Math. Helv. 67 (1992), 664-669.







\bibitem{MaS1}
{\sc P. Mankiewicz and C. Sch\"utt},
{\em A simple proof of an estimate for the approximation
of the Euclidean ball and the Delone triangulation numbers},
 Journal of Approximation Theory 107 (2000), 268-280.


\bibitem{MaS2}
{\sc P. Mankiewicz and C. Sch\"utt},
{\em On the Delone triangulations numbers},
 Journal of Approximation Theory 111 (2001), 139-142.



\bibitem{McVi}
{\sc McClure and R. Vitale},
{\em Polygonal approximation of plane convex bodies}
 J. Math. Anal. Appl. 51 (1975), 326-358.


\bibitem{MePa}
{\sc M. Meyer and A. Pajor},
{\em On the Blaschke-Santal\'o inequality},
Archiv der Mathematik 55 (1990), 82-93.




\bibitem{MeRei2}
{\sc M. Meyer and S. Reisner},
{\em A geometric property of the boundary of symmetric convex bodies
and convexity of floatation surfaces},
Geometriae Dedicata 37 (1991), 327-337.






\bibitem{MiPa}
{\sc V. Milman and A. Pajor},
{\em Isotropic position and inertia ellipsoids and zonoids
of the unit ball of a normed $n$-dimensional space},
Geometric Aspects of Functional Analysis (GAFA)
1987-88, edited by J. Lindenstrauss and V.D. Milman,
Springer-Verlag 1989, 64-104.





\bibitem{Paouris2006} 
{\sc G. Paouris},
{\em Concentration of mass on convex bodies}, Geometric and
Functional Analysis 16,  (2006), 1021--1049.



\bibitem{Pe}
{\sc B. Petkantschin},
{\em Zusammenh\"ange zwischen den Dichten der linearen
Unterr\"aume im $n$-dimensionalen Raume},
 Abhandlungen des Mathematischen Seminars der
Universit\"at Hamburg 11 (1936), 249-310.

\bibitem{Petty}
{\sc C. Petty},
{\em Geominimal surface area},
Geometriae Dedicata 3 (1974), 77-97.

\bibitem{Petty1}
{\sc C. Petty},
{\em Affine isoperimetric problems},
In: Discrete Geometry and Convexity,
J.E. Goodman, E. Lutwak, J. Malkevitch and R. Pollak, eds.,
Ann. New York Acad, Sci 440 (1985), 113-127.








\bibitem{Rei1}
{\sc M. Reitzner},
{\em Random points on the boundary of smooth convex bodies},
 Trans. Amer. Math. Soc., 354 (2002) 2243-2278.


\bibitem{ReSu1} 
{\sc A. R\'enyi and R. Sulanke},
{\em  \"Uber die konvexe H\"ulle von n zuf\"allig 
   gew\"ahlten Punkten}
 Zeitschrift f\"ur Wahrscheinlichkeitstheorie und verwandte Gebiete
2 (1963), 75-84.


\bibitem{ReSu2} 
{\sc  A. R\'enyi and R. Sulanke},
{\em \"Uber die konvexe H\"ulle von n zuf\"allig
   gew\"ahlten Punkten II}
 Zeitschrift f\"ur Wahrscheinlichkeitstheorie und verwandte Gebiete
3  (1964), 138-147



\bibitem{Roc}
{\sc R.T. Rockafellar},
{\em Convex Analysis},
Princeton University Press 1970.



\bibitem{San}
{\sc L.A. Santal\'o},
{\em Integral Geometry and Geometric Probability},
 Addison-Wesley,
 Encyclopedia of Mathematics and its Applications
 1976.











\bibitem{Schneider1}
{\sc R. Schneider},
{\em Affine surface area and convex bodies of elliptic type},
Periodica  Mathematica Hungarica 69 (2014), 120-125.

\bibitem{SchneiderBuch} {\sc R. Schneider}, 
{\em Convex Bodies: The
Brunn-Minkowski Theory}, Encyclopedia of Mathematics and its
Applications  44, Cambridge University Press, Cambridge (1993)



\bibitem{Schue1}
{\sc C. Sch\"utt},
{\em Random polytopes and affine surface area},
 Mathematische Nachrichten 170 (1994), 227-249.





\bibitem{Schue3}
{\sc C. Sch\"utt},
{\em On the affine surface area},
Proceedings of the American Mathematical Society 118
(1993), 1213-1218.


\bibitem{ScWe1}
{\sc C. Sch\"utt and E. Werner},
{\em The convex floating body},
Mathematica Scandinavica 66 (1990), 275-290.




\bibitem{SchueW3}
{\sc C. Sch\"utt and E. Werner},
{\em Random polytopes with vertices on the boundary of a convex body},
Comptes Rendus de l'Acad\'emie des Sciences Paris 331(2000), 697-201.

\bibitem{SchueW4} 
{\sc C. Sch\"utt and E. Werner}, 
{\em Polytopes with vertices chosen randomly from the boundary of a convex body},
Israel  Seminar 2001-2002, Lecture Notes in Mathematics 1807
(V.D. Milman and G. Schechtman, eds.), Springer-Verlag, 2003, 241--422.

\bibitem{SWBook} 
{\sc C. Sch\"utt and E. Werner}, 
{\em Floating Body}, in preparation.







































 \bibitem{Werner1994} 
 {\sc E. Werner}, 
 {\em  Illumination bodies and affine surface area}, Studia
Math. {\bf 110}  (1994), 257-269.


\end{thebibliography}
\end{document}